# A GENERALIZATION OF EULER'S THEOREM ON CONGRUENCIES

Florentin Smarandache
University of New Mexico
200 College Road
Gallup, NM 87301, USA
E-mail: smarand@unm.edu

In the paragraphs which follow we will prove a result which replaces the theorem of Euler:

"If $(a,m)=1$, then $a^{\varphi(m)} \equiv 1(\bmod m)$",

for the case when $a$ and $m$ are not relatively primes.

**A. Introductory concepts.**

One supposes that $m>0$. This assumption will not affect the generalization, because Euler's indicator satisfies the equality:

$\varphi(m)=\varphi(-m)$ (see [1]), and that the congruencies verify the following property:

$a \equiv b(\bmod m) \Leftrightarrow a \equiv b\left(\bmod (-m)\right)$ (see [1] pp 12-13).

In the case of congruence modulo 0, there is the relation of equality. One denotes $(a,b)$ the greatest common factor of the two integers $a$ and $b$, and one chooses $(a,b)>0$.

B. **Lemmas, theorem.**

**Lemma 1:** Let be $a$ an integer and $m$ a natural number $>0$. There exist $d_0, m_0$ from **N** such that $a=a_0d_0$, $m=m_0d_0$ and $(a_0,m_0)=1$.

*Proof:*

It is sufficient to choose $d_0=(a,m)$. In accordance with the definition of the greatest common factor (GCF), the quotients of $a_0$ and $m_0$ and of $a$ and $m$ by their GCF are relatively primes (see [3], pp. 25-26).

**Lemma 2:** With the notations of lemma 1, if $d_0 \neq 1$ and if:

$d_0=d_0^1d_1$, $m_0=m_1d_1$, $(d_0^1,m_1)=1$ and $d_1\neq 1$, then $d_0>d_1$ and $m_0>m_1$, and if $d_0=d_1$, then after a limited number of steps $i$ one has $d_0>d_{i+1}=(d_i,m_i)$.

*Proof:*

$$(0)\begin{cases} a=a_0d_0 \quad ; \quad (a_0,m_0)=1 \\ m=m_0d_0 \quad ; \quad d_0\neq 1 \end{cases}$$

$$(1)\begin{cases} d_0 = d_0^1 d_1 & ; \quad (d_0^1, m_1) = 1 \\ m_0 = m_1 d_1 & ; \quad d_1 \neq 1 \end{cases}$$

From (0) and from (1) it results that $a = a_0 d_0 = a_0 d_0^1 d_1$ therefore $d_0 = d_0^1 d_1$ thus $d_0 > d_1$ if $d_0^1 \neq 1$.

From $m_0 = m_1 d_1$ we deduct that $m_0 > m_1$.
If $d_0 = d_1$ then $m_0 = m_1 d_0 = k \cdot d_0^z$ ($z \in \mathbf{N}^*$ and $d_0 \nmid k$).
Therefore $m_1 = k \cdot d_0^{z-1}$; $d_2 = (d_1, m_1) = (d_0, k \cdot d_0^{z-1})$. After $i = z$ steps, it results $d_{i+1} = (d_0, k) < d_0$.

**Lemma 3:** For each integer $a$ and for each natural number $m > 0$ one can build the following sequence of relations:

$$(0)\begin{cases} a = a_0 d_0 & ; \quad (a_0, m_0) = 1 \\ m = m_0 d_0 & ; \quad d_0 \neq 1 \end{cases}$$

$$(1)\begin{cases} d_0 = d_0^1 d_1 & ; \quad (d_0^1, m_1) = 1 \\ m_0 = m_1 d_1 & ; \quad d_1 \neq 1 \end{cases}$$

……………………………………

$$(s-1)\begin{cases} d_{s-2} = d_{s-2}^1 d_{s-1} & ; \quad (d_{s-2}^1, m_{s-1}) = 1 \\ m_{s-2} = m_{s-1} d_{s-1} & ; \quad d_{s-1} \neq 1 \end{cases}$$

$$(s)\begin{cases} d_{s-1} = d_{s-1}^1 d_s & ; \quad (d_{s-1}^1, m_s) = 1 \\ m_{s-1} = m_s d_s & ; \quad d_s \neq 1 \end{cases}$$

*Proof:*

One can build this sequence by applying lemma 1. The sequence is limited, according to lemma 2, because after $r_1$ steps, one has $d_0 > d_{r_1}$ and $m_0 > m_{r_1}$, and after $r_2$ steps, one has $d_{r_1} > d_{r_1+r_2}$ and $m_{r_1} > m_{r_1+r_2}$, etc., and the $m_i$ are natural numbers. One arrives at $d_s = 1$ because if $d_s \neq 1$ one will construct again a limited number of relations $(s+1), \ldots, (s+r)$ with $d_{s+r} < d_s$.

**Theorem:** Let us have $a, m \in \mathbf{Z}$ and $m \neq 0$. Then $a^{\varphi(m_s)+s} \equiv a^s \pmod{m}$ where s and $m_s$ are the same ones as in the lemmas above.

*Proof:*

Similar with the method followed previously, one can suppose $m > 0$ without reducing the generality. From the sequence of relations from lemma 3, it results that:

(0) (1) (2) (3) (s)

$a = a_0 d_0 = a_0 d_0^1 d_1 = a_0 d_0^1 d_1^1 d_2 = \ldots = a_0 d_0^1 d_1^1 \ldots d_{s-1}^1 d_s$

and

(0) (1) (2) (3) (s)

$m = m_0 d_0 = m_1 d_1 d_0 = m_2 d_2 d_1 d_0 = \ldots = m_s d_s d_{s-1} \ldots d_1 d_0$

and

$m_s d_s d_{s-1} \ldots d_1 d_0 = d_0 d_1 \ldots d_{s-1} d_s m_s$.

From (0) it results that $d_0 = (a, m)$, and from $(i)$ that $d_i = (d_{i-1}, m_{i-1})$, for all $i$ from $\{1, 2, \ldots, s\}$.

$d_0 = d_0^1 d_1^1 d_2^1 \ldots\ldots d_{s-1}^1 d_s$

$d_1 = d_1^1 d_2^1 \ldots\ldots d_{s-1}^1 d_s$

……………………

$d_{s-1} = \qquad d_{s-1}^1 d_s$

$d_s = \qquad d_s$

Therefore

$d_0 d_1 d_2 \ldots\ldots d_{s-1} d_s = (d_0^1)^1 (d_1^1)^2 (d_2^1)^3 \ldots (d_{s-1}^1)^s (d_s^1)^{s+1} = (d_0^1)^1 (d_1^1)^2 (d_2^1)^3 \ldots (d_{s-1}^1)^s$

because $d_s = 1$.

Thus $m = (d_0^1)^1 (d_1^1)^2 (d_2^1)^3 \ldots (d_{s-1}^1)^s \cdot m_s$;

therefore $m_s \mid m$;

(s) (s)

$(d_s, m_s) = (1, m_s)$ and $(d_{s-1}^1, m_s) = 1$

(s-1)

$1 = (d_{s-2}^1, m_{s-1}) = (d_{s-2}^1, m_s d_s)$ therefore $(d_{s-2}^1, m_s) = 1$

(s-2)

$1 = (d_{s-3}^1, m_{s-2}) = (d_{s-3}^1, m_{s-1} d_{s-1}) = (d_{s-3}^1, m_s d_s d_{s-1})$ therefore $(d_{s-3}^1, m_s) = 1$

………..

(i+1)

$1 = (d_i^1, m_{i+1}) = (d_i^1, m_{i+1} d_{i+2}) = (d_i^1, m_{i+3} d_{i+3} d_{i+2}) = \ldots =$

$= (d_i^1, m_s d_s d_{s-1} \ldots d_{i+2})$ thus $(d_i^1, m_s) = 1$, and this is for all $i$ from $\{0, 1, \ldots, s-2\}$.

………..

(0)

$1 = (a_0, m_0) = (a_0, d_1 \ldots d_{s-1} d_s m_s)$ thus $(a_0, m_s) = 1$.

From the Euler's theorem results that:

$(d_i^1)^{\varphi(m_s)} \equiv 1 (\bmod m_s)$ for all $i$ from $\{0, 1, \ldots, s\}$,

$a_0^{\varphi(m_s)} \equiv 1 (\bmod m_s)$

but $a_0^{\varphi(m_s)} = a_0^{\varphi(m_s)}(d_0^1)^{\varphi(m_s)}(d_1^1)^{\varphi(m_s)}...(d_{s-1}^1)^{\varphi(m_s)}$

therefore $a^{\varphi(m_s)} \equiv \underbrace{1........1}_{s+1\ times}(\operatorname{mod} m_s)$

$a^{\varphi(m_s)} \equiv 1(\operatorname{mod} m_s)$.

$a_0^s(d_0^1)^{s-1}(d_1^1)^{s-2}(d_2^1)^{s-3}...(d_{s-2}^1)^1 \cdot a^{\varphi(m_s)} \equiv a_0^s(d_0^1)^{s-1}(d_1^1)^{s-2}...(d_{s-2}^1)^1 \cdot 1(\operatorname{mod} m_s)$.

Multiplying by:

$(d_0^1)^1(d_1^1)^2(d_2^1)^3...(d_{s-2}^1)^{s-1}(d_{s-1}^1)^s$ we obtain:

$a_0^s(d_0^1)^s(d_1^1)^s...(d_{s-2}^1)^s(d_{s-1}^1)^s a^{\varphi(m_s)} \equiv$

$\equiv a_0^s(d_0^1)^s(d_1^1)^s...(d_{s-2}^1)^s(d_{s-1}^1)^s(\operatorname{mod}(d_0^1)^1...(d_{s-1}^1)^s m_s)$

but $a_0^s(d_0^1)^s(d_1^1)^s...(d_{s-1}^1)^s \cdot a^{\varphi(m_s)} = a^{\varphi(m_s)+s}$ and $a_0^s(d_0^1)^s(d_1^1)^s...(d_{s-1}^1)^s = a^s$

therefore $a^{\varphi(m_s)+s} \equiv a^s(\operatorname{mod} m)$, for all $a, m$ from $\mathbf{Z}$ $(m \neq 0)$.

**Observations:**

**(1)** If $(a,m)=1$ then $d=1$. Thus $s=0$, and according to our theorem one has $a^{\varphi(m_0)+0} \equiv a^0(\operatorname{mod} m)$ therefore $a^{\varphi(m_0)+0} \equiv 1(\operatorname{mod} m)$.

But $m = m_0 d_0 = m_0 \cdot 1 = m_0$. Thus:

$a^{\varphi(m)} \equiv 1(\operatorname{mod} m)$, and one obtains Euler's theorem.

**(2)** Let us have $a$ and $m$ two integers, $m \neq 0$ and $(a,m) = d_0 \neq 1$, and $m = m_0 d_0$.

If $(d_0, m_0) = 1$, then $a^{\varphi(m_0)+1} \equiv a(\operatorname{mod} m)$.

Which, in fact, it results from our theorem with $s=1$ and $m_1 = m_0$.

This relation has a similar form to Fermat's theorem:

$$a^{\varphi(p)+1} \equiv a(\operatorname{mod} p).$$

## C. AN ALGORITHM TO SOLVE CONGRUENCIES

One will construct an algorithm and will show the logic diagram allowing to calculate $s$ and $m_s$ of the theorem.

Given as input: two integers $a$ and $m$, $m \neq 0$.

It results as output: $s$ and $m_s$ such that

$a^{\varphi(m_s)+s} \equiv a^s(\operatorname{mod} m)$.

*Method:*

(1) $A := a$
$M := m$
$i := 0$

(2) Calculate $d = (A, M)$ and $M' = M / d$.

(3) If $d = 1$ take $S = i$ and $m_s = M'$ stop.
If $d \neq 1$ take $A := d$, $M = M'$
$i := i + 1$, and go to (2).

Remark: the accuracy of the algorithm results from lemma 3 end from the theorem.
See the flow chart on the following page.
In this flow chart, the SUBROUTINE LCD calculates $D=(A,M)$ and chooses $D>0$.

**Application:** In the resolution of the exercises one uses the theorem and the algorithm to calculate $s$ and $m_s$.

*Example:* $6^{25604} \equiv ?(\text{mod}\,105765)$

One cannot apply Fermat or Euler because (6,105765)=3$\neq 1$. One thus applies the algorithm to calculate $s$ and $m_s$ and then the previous theorem:

$d_0=(6,105765)=3$ $\qquad m_0=105765/3=35255$

$i=0; 3\neq 1$ thus $i=0+1=1$, $d_1=(3,35255)=1$, $m_1=35255/1=35255$.

Therefore $6^{\varphi(35255)+1}\equiv 6^1(\text{mod}\,105765)$ thus $6^{25604}\equiv 6^4(\text{mod}\,105765)$.

## Flow chart:

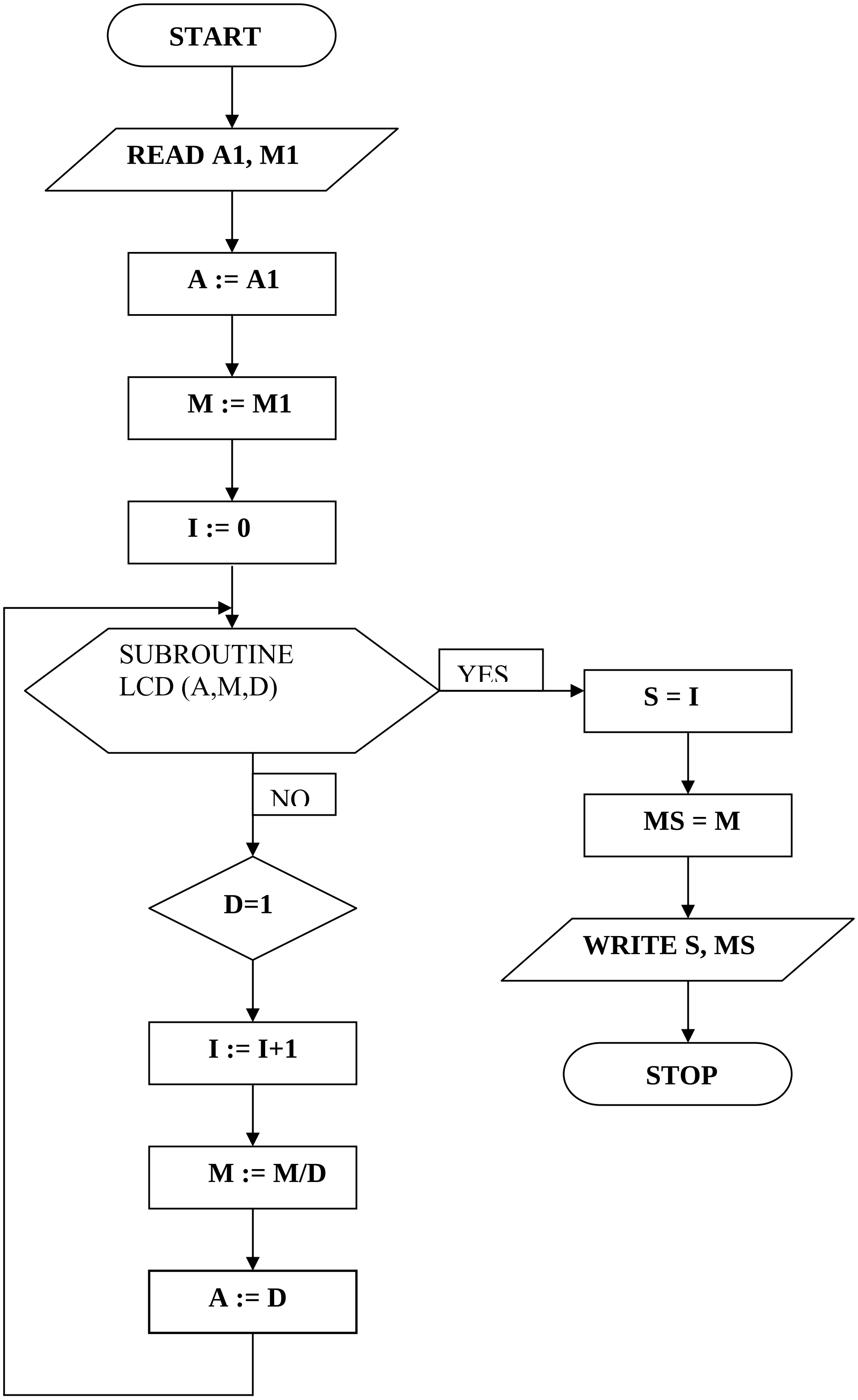